\documentclass[12pt,reqno]{amsart}

\newcommand{\df}{\dfrac}
\newcommand{\tf}{\tfrac}

 \renewcommand{\a}{\alpha}
\renewcommand{\b}{\beta}

\renewcommand{\(}{\left\(}
\renewcommand{\)}{\right\)}
\renewcommand{\[}{\left\[}
\renewcommand{\]}{\right\]}
\renewcommand{\i}{\infty}

\numberwithin{equation}{section}
 \theoremstyle{plain}
\newtheorem{theorem}{Theorem}[section]
\newtheorem{lemma}[theorem]{Lemma}

\begin{document}

\begin{center}

\textbf{A TRANSFORMATION FORMULA INVOLVING THE GAMMA
  AND RIEMANN ZETA FUNCTIONS IN RAMANUJAN'S LOST NOTEBOOK}

\vspace*{.2in}

BRUCE C.~BERNDT\footnotemark[1] AND ATUL DIXIT

\vspace*{.2in}

\centerline{\emph{In Memory of Alladi Ramakrishnan}}

\end{center}

\author{Bruce C.~Berndt}
\address{Department of Mathematics, University of Illinois, 1409 West Green
Street, Urbana, IL 61801, USA} \email{berndt@illinois.edu}
\author{Atul Dixit}
\address{Department of Mathematics, University of Illinois, 1409 West Green
Street, Urbana, IL 61801, USA} \email{aadixit2@illinois.edu}


\footnotetext[1]{Research partially supported by  grant
H98230-07-1-0088.}

\section{Introduction}

Pages 219--227 in the volume \cite{lnb} containing Ramanujan's
lost notebook are devoted to material ``Copied from the Loose
Papers."  We emphasize that these pages are \emph{not} part of the
original lost notebook found by George Andrews at Trinity College
Library, Cambridge in the spring of 1976.  These ``loose papers'',
in the handwriting of G.N.~Watson, are found in the Oxford
University Library, and evidently the original pages in
Ramanujan's handwriting are no longer extant. Most of these nine
pages, which are divided into four rough, partial manuscripts, are
connected with material in Ramanujan's published papers. However,
there is much that is new in these fragments, which will be
completely examined in \cite{geabcbIV}. One claim in the first
manuscript on pages 219--220 is the subject of this short note and
is  the most interesting theorem in the manuscript.  This claim
provides a beautiful series transformation involving the
logarithmic derivative of the gamma function and the Riemann zeta
function.
  To state Ramanujan's claim, it will be convenient to use the familiar notation \cite[p.~952,
formulas 8.360, 8.362, no.~1]{gr}
\begin{equation}\label{w1.15b}
\psi(x):=\df{\Gamma^\prime(x)}{\Gamma(x)}=-\gamma-\sum_{k=0}^{\infty}\left(\df{1}{k+x}-\df{1}{k+1}\right),
\end{equation}
where $\gamma$ denotes Euler's constant.  We also need to recall
the following functions associated with Riemann's zeta function
$\zeta(s)$. Let
\begin{equation*}
\xi(s):=(s-1)\pi^{-\tf{1}{2}s}\Gamma(1+\tf{1}{2}s)\zeta(s).
\end{equation*}
Then Riemann's $\Xi$-function is defined by
\begin{equation*}
\Xi(\tf{1}{2}t):=\xi(\tf{1}{2}+\tf{1}{2}it).
\end{equation*}

\begin{theorem}\label{entry1} Define
\begin{equation}\label{w1.27}
\phi(x):=\psi(x)+\df{1}{2x}-\log x.
\end{equation}
If $\a$ and $\b$ are positive numbers such that $\a\b=1$, then
\begin{multline}\label{w1.26}
\sqrt{\a}\left\{\df{\gamma-\log(2\pi\a)}{2\a}+\sum_{n=1}^{\infty}\phi(n\a)\right\}
=\sqrt{\b}\left\{\df{\gamma-\log(2\pi\b)}{2\b}+\sum_{n=1}^{\infty}\phi(n\b)\right\}\\
=-\df{1}{\pi^{3/2}}\int_0^{\infty}\left|\Xi\left(\df{1}{2}t\right)\Gamma\left(\df{-1+it}{4}\right)\right|^2
\df{\cos\left(\tf{1}{2}t\log\a\right)}{1+t^2}dt,
\end{multline}
where $\gamma$ denotes Euler's constant and $\Xi(x)$ denotes
Riemann's $\Xi$-function.
\end{theorem}

The first identity in \eqref{w1.26} is beautiful in its elegant
symmetry and surprising as well, because why would subtracting the
two leading terms in the asymptotic expansion of the logarithmic
derivative of the Gamma function, in order to gain convergence of
the infinite series on the left side, yield a ``modular relation"
for the resulting function? The second identity in \eqref{w1.26}
is also surprising, for why would the first identity foreshadow a
connection with the Riemann zeta function in the second?

Although Ramanujan does not provide a proof of \eqref{w1.26}, he
does indicate that \eqref{w1.26} ``can be deduced from''
\begin{equation}\label{w1.17}
\int_0^{\infty}\left(\psi(1+x)-\log x\right)\cos(2\pi{nx})dx
=\df{1}{2}\left(\psi(1+n)-\log n\right).
\end{equation}
This latter result was rediscovered by A.P.~Guinand  \cite{apg2}
in 1947, and he later found a simpler proof of this result in
\cite{apg3}. In the footnote at the end of his paper \cite{apg3},
Guinand remarks that T.A.~Brown had told him that he himself had
proved the self-reciprocality of $\psi(1+x)-\log x$ some years
ago, and that when he (Brown) communicated the result to
G.H.~Hardy, Hardy told him that the result was also given by
Ramanujan in a progress report to the University of Madras, but
was not published elsewhere. However, we cannot find this result
in any of the three \emph{Quarterly Reports} that Ramanujan
submitted to the University of Madras \cite{quarterly}, \cite{I}.
Therefore, Hardy's memory was perhaps imperfect; it would appear
that he saw \eqref{w1.17} in the aforementioned manuscript that
Watson had copied.  On the other hand, the only copy of
Ramanujan's \emph{Quarterly Reports} that exists is in Watson's
handwriting! It could be that the manuscript on pages 219--220 of
\cite{lnb}, which is also in Watson's handwriting,  was somehow
separated from the original \emph{Quarterly Reports}, and
therefore that Hardy was indeed correct in his assertion!

The first equality in (\ref{w1.26}) was rediscovered by Guinand in
\cite{apg2} and appears in a footnote on the last page of his
paper \cite[p.~18]{apg2}.  It is interesting that Guinand remarks,
``This formula also seems to have been overlooked."  Here then is
one more instance in which a mathematician thought that his or her
theorem was new, but unbeknownst to him or her, Ramanujan had
beaten them to the punch! We now give Guinand's version of
(\ref{w1.26}).
\begin{theorem} \label{guiver}
For any complex $z$ such that $|arg$ $z|<\pi$, we have
\begin{multline}\label{guiver1}
\sum_{n=1}^{\infty}\left(\frac{\Gamma^{'}}{\Gamma}(nz)-\log nz +\frac{1}{2nz}\right)
+\frac{1}{2z}(\gamma-\log 2\pi z)\\
=\frac{1}{z}\sum_{n=1}^{\infty}\left(\frac{\Gamma^{'}}{\Gamma}\left(\frac{n}{z}\right)
-\log \frac{n}{z}
+\frac{z}{2n}\right)+\frac{1}{2}\left(\gamma-\log
\frac{2\pi}{z}\right).
\end{multline}
\end{theorem}
Although not offering a proof of \eqref{guiver1} in \cite{apg2},
Guinand did remark that it can be obtained by using an appropriate
form of Poisson's summation formula, namely the form given in
Theorem $1$ in \cite{apg1}. Later Guinand gave another proof of
Theorem \ref{guiver} in \cite{apg3}, while also giving extensions
of (\ref{guiver1})
 involving derivatives of the $\psi$-function. He also established a finite version of (\ref{guiver1}) in \cite{apg4}.
 However, Guinand apparently did not discover the connection of his work with Ramanujan's integral
 involving Riemann's $\Xi$-function.

 In this paper we first provide a proof of both identities in Theorem \ref{entry1}.
 In Section $4$, we construct a second proof of (\ref{guiver1}) along the lines suggested by
 Guinand in \cite{apg2}.
We can also provide another proof of (\ref{w1.26}) employing both
\eqref{w1.17} and
\begin{equation}\label{w1.15}
\int_0^{\infty}\left(\df{1}{e^{2\pi{x}}-1}-\df{1}{2\pi{x}}\right)e^{-2\pi{nx}}dx
=\df{1}{2\pi}\left(\log n-\psi(1+n)\right),
\end{equation}
which can be derived from an integral evaluation in \cite[p.~377,
formula 3.427, no.~7]{gr}.
However, this proof is similar but slightly more complicated than
the first proof that we provide below.

Although the Riemann zeta function appears at various instances
throughout Ramanujan's notebooks \cite{nb} and lost notebook
\cite{lnb}, he only wrote one paper in which the zeta function
plays the leading role \cite{riemann}, \cite[pp.~72--77]{cp}.  In
fact, a result proved by Ramanujan in \cite{riemann}, namely
equation (\ref{Ramint}) in Section $3$ below, is a key to proving
\eqref{w1.26}. About the integral involving Riemann's
$\Xi$-function in this result, Hardy \cite{ghh} comments that
``the properties of this integral resemble those of one which
Mr.~Littlewood and I have used, in a paper to be published shortly
in the \textit{Acta Mathematica}, to prove that
\begin{equation}\label{zetamoment}
\int_{-T}^{T}\left|\zeta\left(\frac{1}{2}+ti\right)\right|^2\, dt \sim \frac{2}{\pi} T\log T.\text{''}
\end{equation}

 It is also interesting that on a page in the original lost
notebook \cite[p.~195]{lnb}, Ramanujan defines
 \begin{equation}\label{w1.27a}
\phi(x):=\psi(x)+\df{1}{2x}-\gamma-\log x
\end{equation}
and then concludes that \eqref{w1.26} is valid.  However, with the
definition \eqref{w1.27a} of $\phi(x)$, the series in
\eqref{w1.26} do not converge.  For a more complete discussion of
Ramanujan's incorrect claim, see \cite{geabcbIV}.

\section{Preliminary Results}

We first collect several well-known theorems that we use in our
proof. First, from \cite[p.~191]{con}, for $t\neq 0$,
\begin{equation}\label{cotid}
\sum_{n=1}^{\infty}\frac{1}{t^2+4n^2\pi^2}=\frac{1}{2t}\left(\frac{1}{e^{t}-1}-\frac{1}{t}+\frac{1}{2}\right).
\end{equation}
Second, from \cite[p.~251]{ww}, we find that, for $\text{Re }z>0$,
\begin{equation}\label{w1.28}
\phi(z)=-2\int_{0}^{\infty}\frac{tdt}{(t^2+z^2)(e^{2\pi t}-1)}.
\end{equation}
Third, we require Binet's integral for $\log\Gamma(z)$, i.e., for
$\text{Re } z>0$ \cite[p.~249]{ww}, \cite[p.~377, formula 3.427,
no.~4]{gr},
\begin{equation}\label{w1.30}
\log\Gamma(z)=\left(z-\df{1}{2}\right)\log z -z
+\df{1}{2}\log(2\pi)+\int_0^{\infty}\left(\df{1}{2}-\df{1}{t}+\df{1}{e^t-1}\right)\df{e^{-zt}}{t}dt.
\end{equation}
Fourth, from \cite[p.~377, formula 3.427, no.~2]{gr}, we find that
\begin{equation}\label{w1.29}
\int_{0}^{\infty}\left(\frac{1}{1-e^{-x}}-\frac{1}{x}\right)e^{-x}dx=\gamma,
\end{equation}
where $\gamma$ denotes Euler's constant. Fifth, by Frullani's
integral \cite[p.~378, formula 3.434, no.~2]{gr},
\begin{equation}\label{w1.frullani}
\int_0^{\infty}\df{e^{-\mu{x}}-e^{-\nu{x}}}{x}dx=\log\df{\nu}{\mu},
\qquad \mu, \nu >0.
\end{equation}

\section{First Proof of Theorem \ref{entry1}}

\begin{proof}
Our first goal is to establish an integral representation for the
far left side of \eqref{w1.26}.  Replacing $z$ by $n\a$ in
\eqref{w1.28} and summing on $n$, $1\leq n<\i$, we find, by
absolute convergence, that
\begin{align}\label{phinal}
\sum_{n=1}^{\infty}\phi(n\alpha)&=-2\sum_{n=1}^{\infty}\int_{0}^{\infty}\frac{tdt}{(t^2+n^2\alpha^2)(e^{2\pi
t}-1)}
\notag\\
&=\frac{-2}{\alpha^2}\int_{0}^{\infty}\frac{t}{(e^{2\pi
t}-1)}\sum_{n=1}^{\infty}\frac{1}{({t/\alpha})^2+n^2}.
\end{align}
Invoking \eqref{cotid} in \eqref{phinal}, we see that
\begin{equation}\label{sumphial}
\sum_{n=1}^{\infty}\phi(n\alpha)=-\frac{2\pi}{\alpha}\int_{0}^{\infty}\frac{1}{(e^{2\pi
t}-1)} \left(\frac{1}{e^{2\pi t/\alpha}-1}-\frac{\alpha}{2\pi
t}+\frac{1}{2}\right)dt.
\end{equation}

Next, setting $x=2\pi{t}$ in \eqref{w1.29}, we readily find that
\begin{equation}\label{w1.31}
\gamma=\int_{0}^{\infty}\left(\frac{2\pi}{e^{2\pi
t}-1}-\frac{e^{-2\pi t}}{t}\right)dt.
\end{equation}
By Frullani's integral \eqref{w1.frullani},
\begin{equation}\label{w1.32}
\int_{0}^{\infty}\frac{e^{-t/\alpha}-e^{-2\pi
t}}{t}dt=\log\left(\frac{2\pi}{1/\alpha}\right)=\log(2\pi\alpha).
\end{equation}
Combining \eqref{w1.31} and \eqref{w1.32}, we arrive at
\begin{equation}\label{impint}
\gamma-\log\left(2\pi\alpha\right)=\int_{0}^{\infty}\left(\frac{2\pi}{e^{2\pi
t}-1}-\frac{e^{-t/\alpha}}{t}\right)dt.
\end{equation}

Hence, from \eqref{sumphial} and \eqref{impint}, we deduce that
\begin{align}\label{exp1}
&\sqrt{\alpha}\left(\frac{\gamma -\log(2\pi\alpha)}{2\alpha}+\sum_{n=1}^{\infty}\phi(n\alpha)\right)\\
&\hspace{1cm}=\frac{1}{2\sqrt{\alpha}}\int_{0}^{\infty}\left(\frac{2\pi}{e^{2\pi
t}-1}-\frac{e^{-t/\alpha}}{t}\right)dt\notag\\
&\quad\hspace{1cm}-\frac{2\pi}{\sqrt{\alpha}}\int_{0}^{\infty}\frac{1}{(e^{2\pi
t}-1)}\left(\frac{1}{e^{2\pi t/\alpha}-1}
-\frac{\alpha}{2\pi t}+\frac{1}{2}\right)dt\nonumber\\
&\hspace{1cm}=\int_{0}^{\infty}\left(\frac{\sqrt{\alpha}}{t(e^{2\pi
t}-1)} -\frac{2\pi}{\sqrt{\a}(e^{2\pi t/\alpha}-1)(e^{2\pi
t}-1)}-\frac{e^{-t/\alpha}}{2t\sqrt{\alpha}}\right)dt.\notag
\end{align}

Now from \cite[p.~260, eqn.~(22)]{riemann} or \cite[p.~77]{cp},
for $n$ real,
\begin{align}\label{Ramint}
&\int_{0}^{\infty}\Gamma\left(\frac{-1+it}{4}\right)\Gamma\left(\frac{-1-it}{4}\right)
\left(\Xi\left(\frac{1}{2}t\right)\right)^2\frac{\cos nt}{1+t^2}dt\nonumber\\
&\hspace{1cm}=\int_{0}^{\infty}\left|\Xi\bigg(\frac{1}{2}t\bigg)\Gamma\bigg(\frac{-1+it}{4}\bigg)\right|^{2}
\frac{\cos nt}{1+t^2}dt\nonumber\\
&\hspace{1cm}=\pi^{3/2}\int_{0}^{\infty}\left(\frac{1}{e^{xe^{n}}-1}-\frac{1}{xe^{n}}\right)
\left(\frac{1}{e^{xe^{-n}}-1}-\frac{1}{xe^{-n}}\right)dx.
\end{align}
Letting $n=\tfrac{1}{2}\log\alpha$ and $x=2\pi t/\sqrt{\alpha}$ in
(\ref{Ramint}), we deduce that
\begin{align}\label{befequRamint2}
&-\frac{1}{\pi^{3/2}}\int_{0}^{\infty}\left|\Xi\bigg(\dfrac{1}{2}t\bigg)\Gamma\bigg(\frac{-1+it}{4}\bigg)\right|^{2}
\frac{\cos(\tfrac{1}{2}t\log\alpha)}{1+t^2}dt\\
&\hspace{.1cm}=-\frac{2\pi}{\sqrt{\alpha}}\int_{0}^{\infty}\left(\frac{1}{e^{2\pi
t}-1}-\frac{1}{2\pi t}\right)
\left(\frac{1}{e^{2\pi t/\alpha}-1}-\frac{\alpha}{2\pi t}\right)dt\nonumber\\
&\hspace{.1cm}=\int_{0}^{\infty}\left(\frac{-2\pi/\sqrt{\alpha}}{(e^{2\pi
t/\alpha}-1)(e^{2\pi t}-1)} +\frac{\sqrt{\alpha}}{t(e^{2\pi
t}-1)}+\frac{1}{t\sqrt{\alpha}(e^{2\pi t/\alpha}-1)}
-\frac{\sqrt{\alpha}}{2\pi t^{2}}\right)dt.\notag
\end{align}

Hence, combining (\ref{exp1}) and (\ref{befequRamint2}), in order
to prove that the far left side of \eqref{w1.26} equals the far
right side of \eqref{w1.26}, we see that it suffices to show that
\begin{align}\label{int1zero}
&\int_{0}^{\infty}\left(\frac{1}{t\sqrt{\alpha}(e^{2\pi
t/\alpha}-1)}-\frac{\sqrt{\alpha}}{2\pi t^{2}}
+\frac{e^{-t/\alpha}}{2t\sqrt{\alpha}}\right)dt \notag\\
&\hspace{1cm}=\df{1}{\sqrt{\a}}\int_0^{\infty}\left(\df{1}{u(e^u-1)}-\df{1}{u^2}+\df{e^{-u/(2\pi)}}{2u}\right)du
 = 0,
\end{align}
where we made the change of variable $u=2\pi{t}/\a$.  In fact,
more generally, we show that
\begin{equation}\label{w1.35}
\int_0^{\infty}\left(\df{1}{u(e^u-1)}-\df{1}{u^2}+\df{e^{-ua}}{2u}\right)du
=-\df{1}{2}\log(2\pi{a}),
\end{equation}
so that if we set $a=1/(2\pi)$ in \eqref{w1.35}, we deduce
\eqref{int1zero}.

Consider the integral, for $t>0$,
\begin{align}\label{w1.34}
F(a,t):&=\int_0^{\infty}\left\{\left(\df{1}{e^u-1}-\df{1}{u}+\df{1}{2}\right)\df{e^{-tu}}{u}
+\df{e^{-ua}-e^{-tu}}{2u}\right\} du\notag\\
&=\log\Gamma(t)-\left(t-\df{1}{2}\right)\log{t}+t-\df{1}{2}\log(2\pi)+\df{1}{2}\log\df{t}{a},
\end{align}
where we applied \eqref{w1.30} and \eqref{w1.frullani}.  Upon the
integration of \eqref{w1.15b}, it is easily gleaned that, as
$t\to0$,
\begin{equation*}
\log\Gamma(t)\sim-\log t-\gamma{t},
\end{equation*}
where $\gamma$ denotes Euler's constant.  Using this in
\eqref{w1.34}, we find, upon simplification, that, as $t\to0$,
\begin{align*}
F(a,t)\sim-\gamma{t}-t\log{t}+t-\df{1}{2}\log(2\pi)-\df{1}{2}\log{a}.
\end{align*}
Hence,
\begin{equation}\label{w1.36}
\lim_{t\to0}F(a,t)=-\df{1}{2}\log(2\pi{a}).
\end{equation}
Letting $t$ approach 0 in \eqref{w1.34}, taking the limit under
the integral sign on the right-hand side using Lebesgue's
dominated convergence theorem, and employing \eqref{w1.36}, we
immediately deduce \eqref{w1.35}.  As previously discussed, this
is sufficient to prove the equality of the first and third
expressions in \eqref{w1.26}, namely,
\begin{equation}\label{w1.good}
\sqrt{\a}\left\{\df{\gamma-\log(2\pi\a)}{2\a}+\sum_{n=1}^{\infty}\phi(n\a)\right\}
=-\df{1}{\pi^{3/2}}\int_0^{\infty}\left|\Xi\left(\df{1}{2}t\right)\Gamma\left(\df{-1+it}{4}\right)\right|^2
\df{\cos\left(\tf{1}{2}t\log\a\right)}{1+t^2}dt.
\end{equation}

Lastly, using \eqref{w1.good} with $\a$ replaced by $\b$ and
employing the relation $\a\b=1$, we conclude that
\begin{align*}
&\sqrt{\b}\left\{\df{\gamma-\log(2\pi\b)}{2\b}+\sum_{n=1}^{\infty}\phi(n\b)\right\}\\
&\hspace{1cm}=-\df{1}{\pi^{3/2}}\int_0^{\infty}\left|\Xi\left(\df{1}{2}t\right)\Gamma\left(\df{-1+it}{4}\right)\right|^2
\df{\cos\left(\tf{1}{2}t\log\b\right)}{1+t^2}dt\\
&\hspace{1cm}=-\df{1}{\pi^{3/2}}\int_0^{\infty}\left|\Xi\left(\df{1}{2}t\right)\Gamma\left(\df{-1+it}{4}\right)\right|^2
\df{\cos\left(\tf{1}{2}t\log(1/\a)\right)}{1+t^2}dt\\
&\hspace{1cm}=-\df{1}{\pi^{3/2}}\int_0^{\infty}\left|\Xi\left(\df{1}{2}t\right)\Gamma\left(\df{-1+it}{4}\right)\right|^2
\df{\cos\left(\tf{1}{2}t\log\a\right)}{1+t^2}dt.
\end{align*}
Hence, the equality of the second and third expressions in
\eqref{w1.26} has been demonstrated, and so the proof is complete.
\end{proof}

\section{Second Proof of \eqref{w1.26}}

In this section we give our second proof of the first identity in
\eqref{w1.26} using Guinand's generalization of Poisson's
summation formula in \cite{apg1}. We emphasize that this route
does not take us to the integral involving Riemann's
$\Xi$-function in the second identity of \eqref{w1.26}. First, we
reproduce the needed version of the Poisson summation formula from
Theorem 1 in \cite{apg1}.

\begin{theorem}\label{genpoi}
If $f(x)$ is an integral, tends to zero at infinity, and
$xf^{'}(x)$ belongs to $L^{p}(0,\infty)$, $(1<p\leq 2)$, then
\begin{equation}
\lim_{N\to\infty}\left(\sum_{n=1}^{N}f(n)-\int_{0}^{N}f(t) dt\right)
=\lim_{N\to\infty}\left(\sum_{n=1}^{N}g(n)-\int_{0}^{N}g(t) dt\right),
\end{equation}
where
\begin{equation}\label{gexform}
g(x)=2\int_{0}^{\to\infty}f(t)\cos(2\pi xt) dt.
\end{equation}
\end{theorem}
Next, we state a lemma\footnote{The authors are indebted to
M.~L.~Glasser for the proof of this lemma. The authors' original
proof of this lemma was substantially longer than Glasser's given
here.} that will subsequently be used in our  proof of
\eqref{w1.26}.

\begin{lemma}\label{caid}
\begin{equation}
\int_{0}^{\infty}\left(\psi(t+1)-\frac{1}{2(t+1)}-\log t\right)dt = \frac{1}{2}\log 2\pi.
\end{equation}
\end{lemma}

\begin{proof}
Let $I$ denote the integral on the left-hand side. Then,
{\allowdisplaybreaks\begin{align}\label{Ieval}
I&=\int_{0}^{\infty}\frac{d}{dt}\left(\ln\frac{e^t\Gamma(t+1)}{t^t\sqrt{t+1}}\right)\, dt\nonumber\\
&=\lim_{t\to\infty}\ln\frac{e^t\Gamma(t+1)}{t^t\sqrt{t+1}}-\lim_{t\to 0}\ln\frac{e^t\Gamma(t+1)}{t^t\sqrt{t+1}}\nonumber\\
&=\ln\lim_{t\to\infty}\frac{e^t\Gamma(t+1)}{t^t\sqrt{t+1}}-\ln\left(\lim_{t\to 0}e^t\Gamma(t+1)\right)
-\lim_{t\to 0}t\ln t-\lim_{t\to 0}\frac{1}{2}\ln(t+1)\nonumber\\
&=\ln\lim_{t\to\infty}\frac{e^t\Gamma(t+1)}{t^t\sqrt{t+1}}.
\end{align}}
Next, Stirling's formula \cite[p.~945, formula 8.327]{gr} tells us that,
\begin{equation}\label{stirling}
\Gamma(z)\sim z^{z-1/2}e^{-z}\sqrt{2\pi},
\end{equation}
as $|z|\to\infty$ and $|\arg$ $z|\leq \pi-\delta$, where $0<\delta<\pi$.
Hence, employing (\ref{stirling}), we find that
\begin{equation}
\frac{e^t\Gamma(t+1)}{t^t\sqrt{t+1}}\sim \left(1+\frac{1}{t}\right)^t\frac{\sqrt{2\pi}}{e},
\end{equation}
so that
\begin{equation}\label{lim}
\lim_{t\to\infty}\frac{e^t\Gamma(t+1)}{t^t\sqrt{t+1}}=\sqrt{2\pi}.
\end{equation}
Thus from (\ref{Ieval}) and (\ref{lim}), we conclude that
\begin{equation}
I=\frac{1}{2}\ln 2\pi.
\end{equation}
\end{proof}

Now we are ready to give our second proof of \eqref{w1.26}.  We
first prove it for real $z>0$. Let
\begin{equation}\label{fex}
f(x)=\psi(xz+1)-\log xz .
\end{equation}
We show that $f(x)$ satisfies the hypotheses of Theorem \ref{genpoi}. From (\ref{w1.17}),
it is true that $f(x)$ is an integral. Next, we need two formulas for $\psi(x)$.

First, from \cite[p.~259, formula 6.3.18]{as} for $|\arg$ $z|<\pi$, as $z\to\infty$,
\begin{equation}\label{psiasymp}
\psi(z) \sim \log z-\frac{1}{2z}-\frac{1}{12z^{2}}+\frac{1}{120z^{4}}-\frac{1}{252z^{6}}+\cdots.
\end{equation}

Second, we use the fact from \cite[p.~250]{ww} that
\begin{equation}\label{psiprime}
\psi^{'}(z)=\sum_{n=0}^{\infty}\frac{1}{(z+n)^2}.
\end{equation}
From (\ref{psiasymp}), we have
\begin{equation}
f(x)\sim\frac{1}{2xz}-\frac{1}{12x^{2}z^{2}}+\frac{1}{120x^{4}z^{4}}-\frac{1}{252x^{6}z^{6}}+\cdots,
\end{equation}
so that
\begin{align}
\lim_{x\to\infty}f(x)&=0.
\end{align}

Next, we show that $x f^{'}(x)$ belongs to $L^{p}(0,\infty)$ for some $p$ such that $1<p\leq 2$.

Using (\ref{psiasymp}), we find that
\begin{align}
xf^{'}(x)=xz\psi^{'}(xz)-\frac{1}{xz}-1 \sim -\frac{1}{2xz},
\end{align}
so that $|xf^{'}(x)|^{p}\sim\tfrac{1}{2^{p}z^{p}x^{p}}$. Now $p>1$ implies that $x f^{'}(x)$ is
locally integrable near $\infty$. Also, using (\ref{psiprime}), we have
\begin{align}
\lim_{x\to 0}xf^{'}(x)&=\lim_{x\to 0}\left(xz\sum_{n=0}^{\infty}\frac{1}{(xz+n)^2}-\frac{1}{xz}-1\right)\nonumber\\
&=\lim_{x\to 0}\left(xz\sum_{n=1}^{\infty}\frac{1}{(xz+n)^2}-1\right)\nonumber\\
&=-1.
\end{align}%
This proves that $x f^{'}(x)$ is locally integrable near $0$. Hence it is true that $x f^{'}(x)$
belongs to $L^{p}(0,\infty)$ for some $p$ such that $1<p\leq 2$.

Now from (\ref{gexform}) and (\ref{fex}), we find that
\begin{equation*}
g(x)=2\int_{0}^{\infty}(\psi(tz+1)-\log tz)\cos\left(2\pi xt\right)dt.
\end{equation*}
Employing a change of variable $y=tz$ and using (\ref{w1.17}), we
find that
\begin{align}\label{gex}
g(x)&=\frac{2}{z}\int_{0}^{\infty}(\psi(y+1)-\log y)\cos\left(2\pi xy/z\right)dy\nonumber\\
&=\frac{1}{z}\left(\psi\left(\frac{x}{z}+1\right)-\log\left(\frac{x}{z}\right)\right).
\end{align}%

Substituting the expressions for $f(x)$ and $g(x)$ from
(\ref{fex}) and (\ref{gex}), respectively, in (\ref{genpoi}), we
find that
\begin{align}
&\lim_{N\to\infty}\left(\sum_{n=1}^{N}(\psi(nz+1)-\log nz)-\int_{0}^{N}(\psi(tz+1)-\log tz) dt\right)\nonumber\\
&=\frac{1}{z}\left[\lim_{N\to\infty}\left(\sum_{n=1}^{N}\left(\psi\left(\frac{n}{z}+1\right)
-\log \frac{n}{z}\right)-\int_{0}^{N}\left(\psi\left(\frac{t}{z}+1\right)-\log \frac{t}{z}\right) dt\right)\right].
\end{align}
Thus,
\begin{align}\label{beffin}
&\lim_{N\to\infty}\left(\sum_{n=1}^{N}\left(\frac{\Gamma^{'}}{\Gamma}(nz)+\frac{1}{2nz}-\log nz\right)
+\sum_{n=1}^{N}\frac{1}{2nz}-\int_{0}^{N}(\psi(tz+1)-\log tz) dt\right)\nonumber\\
&=\frac{1}{z}\left[\lim_{N\to\infty}\left(\sum_{n=1}^{N}\left(\frac{\Gamma^{'}}{\Gamma}\left(\frac{n}{z}\right)
+\frac{z}{2n}-\log \frac{n}{z}\right)+\sum_{n=1}^{N}\frac{z}{2n}-\int_{0}^{N}\left(\psi\left(\frac{t}{z}+1\right)
-\log \frac{t}{z}\right) dt\right)\right].
\end{align}

Now if we can show that
\begin{equation}\label{rem1}
\lim_{N\to\infty}\left(\sum_{n=1}^{N}\frac{1}{2nz}-\int_{0}^{N}(\psi(tz+1)-\log tz) dt\right)
=\frac{\gamma-\log 2\pi z}{2z},
\end{equation}
then replacing $z$ by $1/z$ in (\ref{rem1}) will give us
\begin{equation}\label{rem2}
\lim_{N\to\infty}\left(\sum_{n=1}^{N}\frac{z}{2n}-\int_{0}^{N}\left(\psi\left(\frac{t}{z}+1\right)
-\log \frac{t}{z}\right) dt\right)=\frac{z(\gamma-\log (2\pi/z))}{2}.
\end{equation}
Then substituting (\ref{rem1}) and (\ref{rem2}) in (\ref{beffin}) will complete the proof of the theorem.
 So our goal now is to prove (\ref{rem1}).
{\allowdisplaybreaks\begin{align}
&\lim_{N\to\infty}\left(\sum_{n=1}^{N}\frac{1}{2nz}-\int_{0}^{N}(\psi(tz+1)-\log tz) dt\right)\nonumber\\
&=\lim_{N\to\infty}\left(\frac{1}{2z}\left(\sum_{n=1}^{N}\frac{1}{n}-\log N\right)+\frac{\log N}{2z}
-\int_{0}^{N}(\psi(tz+1)-\log tz) dt\right)\nonumber\\
&=\frac{\gamma}{2z}+\lim_{N\to\infty}\left(-\frac{\log z}{2z}
+\frac{\log Nz}{2z}-\int_{0}^{N}(\psi(tz+1)-\log tz) dt\right)\nonumber\\
&=\frac{\gamma}{2z}-\frac{\log z}{2z}+\lim_{N\to\infty}\left(\frac{\log (Nz+1)}{2z}
-\frac{1}{z}\int_{0}^{Nz}(\psi(t+1)-\log t) dt-\frac{1}{2z}\log\left(1+\frac{1}{Nz}\right)\right)\nonumber\\
&=\frac{\gamma}{2z}-\frac{\log z}{2z}+\frac{1}{z}\lim_{N\to\infty}\left(\frac{\log (Nz+1)}{2}
-\int_{0}^{Nz}(\psi(t+1)-\log t) dt\right)\nonumber\\
&=\frac{\gamma}{2z}-\frac{\log z}{2z}+\frac{1}{z}\lim_{N\to\infty}\left(\frac{1}{2}\int_{0}^{Nz}\frac{1}{t+1} dt
-\int_{0}^{Nz}(\psi(t+1)-\log t) dt\right)\nonumber\\
&=\frac{\gamma}{2z}-\frac{\log z}{2z}-\frac{1}{z}\lim_{N\to\infty}\int_{0}^{Nz}\left(\psi(t+1)-\frac{1}{2(t+1)}
-\log t\right)dt\nonumber\\
&=\frac{\gamma}{2z}-\frac{\log z}{2z}-\frac{1}{z}\int_{0}^{\infty}\left(\psi(t+1)-\frac{1}{2(t+1)}
-\log t\right)dt\nonumber\\
&=\frac{\gamma}{2z}-\frac{\log z}{2z}-\frac{\log 2\pi}{2z}\nonumber\\
&=\frac{\gamma-\log 2\pi z}{2z},
\end{align}}%
where in the antepenultimate line, we have made use of Lemma \ref{caid}.
This completes the proof of (\ref{rem1}) and hence the proof of Theorem \ref{guiver} for real $z>0$.
But both sides of (\ref{guiver1}) are analytic for $|\arg$ $z|<\pi$. Hence by analytic continuation,
 the theorem is true for all complex $z$ such that $|\arg$ $z|<\pi$.

\end{document}